\documentclass[final]{amsart}

\usepackage{xypic}

\newcommand{\Z}{\ensuremath{{\mathbb Z}}}
\newcommand{\Q}{\ensuremath{\mathbb Q}}
\newcommand{\C}{\ensuremath{\mathbb C}}

\newcommand{\coh}[2]{{\rm H}^{#1}(#2)}
\renewcommand{\div}{{\rm Div}}

\newcommand{\isoarrow}{\tilde{\longrightarrow}}
\newcommand{\gr}[2]{{\rm Gr}^{#1}_{#2}}

\newcommand{\struct}[1]{{\mathcal O}_{#1}}
\newcommand{\structd}{{\struct{X}(D)^r}}

\newcommand{\spec}{{\rm Spec}}
\newcommand{\dual}{{\hspace{2pt}\check{}}}
\renewcommand{\hom}[2]{{\rm Hom}_{#1}(#2)}
\newcommand{\ext}[3]{\ensuremath{\rm Ext}^{#1}_{#2}(#3)}

\newtheorem{proposition}{Propositon}[section]
\newtheorem{coro}[proposition]{Corollary}
\newtheorem{prop}[proposition]{Proposition}
\newtheorem{corollary}[proposition]{Corollary}
\newtheorem{theorem}[proposition]{Theorem}
\newtheorem{lemma}[proposition]{Lemma}

\newcommand{\bdiv}{ {\bf Div}}

\newcommand{\E}{\mathcal E}
\newcommand{\F}{\mathcal F}
\newcommand{\G}{\mathcal G}

\newcommand{\bu}{_{\bullet}}
\newcommand{\ubu}{^{\bullet}}

\newcommand{\cosk}{\ensuremath{\rm cosk}}
\newcommand{\sk}{\ensuremath{\rm sk}}

\newcommand{\sX}{{\ensuremath{\mathfrak X}}}
\newcommand{\sY}{{\ensuremath{\mathfrak Y }}}
\newcommand{\sS}{{\ensuremath{\mathfrak S}}}
\newcommand{\mo}{{\mathfrak M}^s}
\newcommand{\bun}{{\rm Bun}}
\newcommand{\utop}{^{\rm top}}
\newcommand{\hilb}{{\rm Hilb}}

\newcommand{\schemes}{{\bf schemes}}
\newcommand{\mtop}{{\bf top}}

\newcommand{\ri}{{\rightarrow}}
\newcommand{\gl}{{\rm GL}}
\newcommand{\pgl}{{\rm PGL}}
\newcommand{\gm}{{\mathbb G}_m}

\newcommand{\ur}{{\underline{r}}}
\newcommand{\ud}{{\underline{d}}}
\newcommand{\flag}{{\rm Flag}}
\newcommand{\prd}{{(\ur,\ud)}}
\newcommand{\divf}{{\rm Div}}
\newcommand{\tE}{\tilde{\E}}

\begin{document}

\title{On the Cohomology of Moduli  of Vector Bundles}

\author[A. Dhillon]{Ajneet Dhillon}

\address{Department of Mathematics \\
        University of British Columbia}
\email{adhillon@math.ubc.ca}

\subjclass{14}

\begin{abstract}
We compute some Hodge and Betti numbers of the moduli space of stable  rank $r$
degree $d$ vector bundles on a smooth projective curve. We
do not assume  $r$ and $d$ are coprime.
In the process we equip the cohomology of an arbitrary algebraic stack with a
functorial mixed Hodge structure. This Hodge structure is 
computed in the case of the moduli stack of rank $r$, degree
$d$ vector bundles on a curve. Our methods also yield a formula 
for the Poincare
polynomial of the moduli stack that is valid over any 
ground field.
\end{abstract}

\maketitle

\section{Introduction}

We will work over a ground field $k$. For an algebraic 
stack $\sY$ defined over $k$, when we speak of its cohomology
will mean its $\ell$-adic cohomology in the smooth topology, except
when $k=\C$, in which case we will mean the cohomology of the 
constant sheaf with values in $\Q$ in the usual topology. These
constructions are reviewed in \S 2. The generic notation
$\coh{*}{\sY}$ for these cohomology theories is used and it will be clear
from the context what is meant. As we are working over
a possibly non algebraically closed field, we remind the
reader that the $\ell$-adic cohomology is always
defined by first passing to an algebraic closure, that is
\[
 {\rm H}_{\rm sm}^i(\sY,\Q_{\ell})  \stackrel{\text{defn}}{=}
 {\rm H}_{\rm sm}^i(\sY\otimes_k \bar{k},\Q_{\ell}).
\]
The ground field $k$ is detected only in the Galois action
on these cohomology groups.

Let $X$ be a smooth, geometrically connected, projective curve, 
defined over $k$, of
genus $g\ge 2$. Fix integers
$r>0$ and $d$ and let $\mo_{r,d}$ be the moduli space of rank $r$
and degree $d$ stable vector bundles on this curve. We denote 
by $\bun_{r,d}$ be the moduli stack of rank $r$ and degree $d$ vector bundles 
on $X$. The integers $r$ and $d$ will frequently be omitted from the
notation. 

We will calculate the Betti numbers, $\dim \coh{i}{\mo}$,
(and Hodge numbers for $k=\C$) when 
$i<2(r-1)(g-1)$ in this article. For $r$ and $d$ coprime this
question has been extensively studied, see \cite{atiyah:82},
\cite{harder:75} and \cite{bifet:94}. On the other hand when
$r$ and $d$ are no longer coprime the question has remained 
open and only partial results exist, which we now describe. 
In rank two, a desingularization $\tilde{{\mathfrak M}^{ss}}$ 
of ${\mathfrak M}^{ss}$ has been constructed by C. Seshadri. In
\cite{balaji:90},\cite{balaji:93} and \cite{balaji:97} its
cohomology is studied. In \cite{arapura:01} the Hodge and
Betti numbers of $\coh{i}{\mo}$ are computed for
\[
 i<2(r-1)g-(r-1)(r^2+3r+1)-7.
\]

Our method is to continue the study of the ind scheme
${\bf \div}$ that was started in \cite{bifet:94}.(See \S 3
for the definition of ${\bf \div}$.)
 In
this cited paper the Poincare polynomial of this ind scheme and
its Shatz strata are computed. We will review this computation
is section 3. 
In section 4 we show that the natural map 
\begin{eqnarray}
\label{E:isom}
 {\bf \bdiv}{\rightarrow}\bun
\end{eqnarray}
is a quasi-isomorphism. This allows us to compute 
the Betti numbers of the stack. 
Over $\C$ this was first done in \cite{atiyah:82}. In this paper
the Poincare polynomial of the classifying space of the 
gauge group is written down. A simple argument shows that
in fact $\bun$ and this classifying space have the same
cohomology. In the introduction to \cite{bifet:94}, the remark
was made that ${\bf \div}$ and this classifying space have 
the same Poincare polynomial and hence this coincidence
is explained by the above isomorphism.

To obtain the Betti numbers of $\mo$ we 
prove a comparison theorem between the cohomology of
$\bun^s$ and $\mo$, see 5. As (\ref{E:isom}) holds for stable
loci, this theorem reduces the
study of the cohomology of the coarse moduli space $\mo$ to that
of the fine moduli space $\bdiv^s$, where superscript $s$ refers
to the stable locus. We are unable to completely describe the
cohomology of this ind scheme so we instead provide
an upper bound on the codimension of the complement of
${\bf Div}^s$ in ${\bf Div}$. 

Although not completely necessary here, it is desirable to
provide
 a suitable theory of mixed Hodge 
structures for algebraic stacks. Our first task will be to sketch 
such a construction.
Note that such a construction was first suggested in
\cite{teleman:98} but has not been published, so it
is provided here.

The construction of a functorial mixed Hodge
structure on the cohomology of a stack
is entirely analogous to that
given \cite{HodgeIII}. Given an algebraic stack
${\mathfrak X}$ and a smooth presentation
\[
 P{\rightarrow}{\mathfrak X},
\]
we can form the simplicial algebraic space whose 
$n$th term is
\[
 \underbrace{P\times_{\mathfrak X}P\times_{\mathfrak X}\ldots
P\times_{\mathfrak X}P}_{n\text{ times}}, 
\]
or in the notation of \cite{HodgeIII}
\[
 {\rm cosk}(P/{\mathfrak X}).
\]
Essentially, the method for equipping such an
algebraic space with a functorial mixed Hodge structure
is given in \cite{HodgeIII}, provided it is separated
and of finite type.
These last two conditions are a hindrance here as the stack 
$\bun$ is neither separated nor of finite type.
To remove the condition of separatedness we construct a
new simplicial scheme $Y\bu$ and a map
\[
 Y\bu{\rightarrow}{\rm cosk}(P/{\mathfrak X}),
\] 
such that $Y\bu$ is separated and the map is of cohomological
descent. The finite type assumption is not really essential
in \cite{HodgeIII}, what is important is that the
cohomology of the stack be finite dimensional.

{\bf{Acknowledgments}} This paper is largely based 
on my PhD thesis. I would like 
to thank Donu Arapura, my advisor, for many 
helpful discussions through the years. This work
would not have possible without the help of 
Pramath Sastry. It was his course that introduced 
me to these ideas. I was possible to simplify some
proofs that appear here thanks to conversations
with Kai Behrend.

\section{Hodge Theory for Algebraic Stacks}

It is not practical to redo the entire contents 
of \cite{HodgeII} and \cite{HodgeIII} here as 
the modifications are only minor. We will therefore
refer to these works for the bulk of the construction.

Firstly we need to recall the definition of  the cohomology of an algebraic
stack ${\mathfrak X}{\rightarrow}{\spec (k)}.$ The stack
$\sX$ is a category fibered over $\schemes /k$. This second
category has a smooth topology, so we define an arrow 
to be a cover if its image in $\schemes /k$ is. This allows
us to consider the $\ell$-adic cohomology in the smooth 
topology on $\sX$. For details, see \cite{behrend:03} or
\cite{laumon:00}.
When
$k=\C$, we may pass to the underlying topological stack
\[
 {\sX}^{\rm top}{\rightarrow}{\mtop}.
\]
Similarly, one may define a Grothendieck topology on this 
stack by use of the big site on
$\mtop$. Given a coefficient ring $F$ we denote by
\[
 \coh{*}{\sX,F},
\]
the cohomology of the constant sheaf with 
values in $F$ on this site. (We remind the reader
of our conventions for when $F=\Q$ stated at the 
beginning of the article). A good introduction
to the cohomology of stacks can be found in 
Kai Behrend's January 2002
talk at MSRI. It
is available on streaming video at
\[
 \text{http://www.msri.org}
\]

These definitions are not completely necessary here, as we 
will be replacing our stack by a simplicial space and
the cohomology of this simplicial space will be the same 
as that of the stack. 

General references for simplicial objects and
cohomological descent are \cite{saint-donat:72} and
\cite{HodgeIII}. 
For a simplicial object, denote by $\sk_n$ the
$n$th truncation functor and by $\cosk_n$ its 
right adjoint. Fix a locally finite stack $\sX$ over
$k$ 
and a smooth presentation
\[
 {\alpha}:P{\rightarrow}\sX.
\]

\begin{proposition} (i)The map ${\alpha}$ is
of universal cohomological descent for the
smooth topology.

\noindent
(ii)
The map
\[
 {\alpha}\utop:P\utop{\rightarrow}\sX\utop
\]
is of universal cohomological descent for
the usual topology.
\end{proposition}

\begin{proof}
The proof of (i) can be found in \cite{behrend:03}.
We give a sketch only of (ii) and leave the details to the 
reader.  Recall that $\sX\utop=[P\utop\bu]$, where
$P\utop\bu$ is the topological space in groupoids
defined by
\[
\xymatrix{
 (P\times_{\sX}P)\utop \ar@<1ex>[r]^(.6)s \ar@<-1ex>[r]_(.6)t & 
 P\utop. }
\]
Both the arrows $s$ and $t$ admit sections locally on $X\utop$ 
as they maps underlying smooth morphisms of algebraic spaces.
Using this fact, one shows that for every topological space $T$ and
every $T{\rightarrow}\sX\utop$ the map
\[
 T\times_{\sX\utop}P\utop{\rightarrow}T
\]
admits sections locally on $T$. Now the result follows 
as the question is local on the base $\sX\utop$.
\end{proof}

\begin{corollary}
The natural augmentation map
\[
 \cosk(P/\sX){\rightarrow}\sX
\]
induces an isomorphism
\[
 \coh{i}{\cosk(P/\sX)}\isoarrow\coh{i}{\sX}.
\]
\end{corollary}

It is worth noting that the following spectral 
sequence relates the cohomology of the components
of $\cosk(P/\sX)$ to that of $\sX$.

\begin{proposition}
\label{P:ss}
Let $Z\bu$ be a simplicial space. Then the there
is a spectral sequence with
\[
 E_1^{pq}=\coh{q}{Z_p}
\]
abutting to $\coh{p+q}{Z\bu}$.
\end{proposition}

\begin{proof}
See \cite{saint-donat:72}
\end{proof}

For the remainder of this section we will take
$k=\C$. 
Let ${\bf lfschemes}/\C$ be the full subcategory 
of $\schemes/\C$ consisting of schemes that
are separated and are disjoint unions of 
schemes of finite type over $\C$. Let
${\bf lfss}_k$ be the category of 
$k$-truncated simplicial objects in
${\bf lfschemes}/\C$. Our next task is to
construct a smooth simplicial scheme $Y\bu$ 
in ${\bf lfss}_{\infty}$ with a map
\[
 Y\bu{\rightarrow}\cosk(P/\sX)
\]
that is a hypercover. First let us recall
the standard method for construction 
of hypercovers.

In what follows, a simplicial space could mean
simplicial scheme, simplicial algebraic space
or a simplicial topological space.

Consider a $m$-truncated simplicial space
$X\bu$ augmented towards a stack $\sS$, i.e
\[
a:X\bu{\rightarrow}\sS.
\]
Recall that $a$ is called a \emph{hypercover} if
the canonical maps deduced from adjunction
\[
 X_{n+1}{\rightarrow}(\cosk\ \sk X\bu)_{n+1}
\quad  \text{for }-1\le n\le m-1,
\]
are of universal cohomological descent. This definition
makes sense for $m=\infty$.
Recall the following (\cite[3.3.3]{saint-donat:72}):

\begin{theorem}
\label{T:hypercover}
If $a:X\bu{\rightarrow}\sS$ is a hypercover as above, then
the natural map
\[
 \cosk(X\bu/\sS){\rightarrow}\sS
\]
is of universal cohomological descent.
\end{theorem}

We describe below the main method for 
constructing hypercovers.

A $k$-truncated simplicial space $X\bu$ is said
to be \emph{split} if there exists for each $j$, $k\ge j\ge 0$,
a subobject $NX_j$ of $X_j$ such that the morphisms

\[
\coprod s:\coprod_{i\le n}
\coprod_{s\in\hom{}{\Delta_n,\Delta_i}} N(X_i) \ri X_n
\]

are isomorphisms, for $n\le k$. This definition makes sense
for $k=\infty$.

Let $X\bu$ be a split $k$-truncated simplicial space with $k$
a finite number. We denote by ${\alpha}(X\bu)$ the 
triple $(X',N,{\beta})$, where

\noindent
(i) the $(k-1)$-truncated simplicial space $X'$ obtained by restricting
$X\bu$.

\noindent
(ii) $N=NX_k$

\noindent
(iii) The canonical map
\[
 {\beta}:NX_k\ri ({\rm cosk}_{k-1}{\rm sk}_{k-1}(X\bu))_k.
\]

The triple ${\alpha}(X)=(X',N,{\beta})$ satisfies
the following condition
\[
\begin{array}{c}
\text{(S)} \quad
X'\text{ is a }(k-1)\text{-truncated split simplicial space and }
{\beta}\text{ is a map } \\
{\beta}:N\ri ({\rm cosk}_{k-1} X')_k
\end{array}
\]

\begin{prop}
\label{P:hyperconstruct}
(i) Let $(X',N,{\beta})$ be a triple satisfying
(S). Upto isomorphism, there exists a unique split 
$k$-truncated $X\bu$ with
\[
{\alpha}(X)\cong(X',N,{\beta})
\] 

\noindent
(ii) In the setup of the previous part suppose
$Z$ is a $k$-truncated simplicial space. To
give a map $f:X{\rightarrow}Z$ is the same as giving 
the following data:

\noindent
(a) a map $f':X'{\rightarrow}\sk_{k-1}(Z)$

\noindent
(b) a map ${f'}':N{\rightarrow}Z_k$ such that the
following diagram commutes:
\[
\xymatrix{ 
N \ar[r] \ar[d] & (\cosk X')_k \ar[d] \\
Z_k \ar[r] & (\cosk \sk_{k-1} Z)_k.
}
\]

\end{prop}

\begin{proof}
This is Proposition 5.1.3 of \cite{saint-donat:72}.
\end{proof}

Now recall our setup, from earlier in this section,
we had a stack $\sX$ and a smooth presentation
\[
 P{\rightarrow}\sX.
\]
We construct our hypercover $Y\bu$ of
$\cosk(P/\sX)$ inductively as follows:

$k=0:$ Let $P\ri\sX$ be a  presentation. We may assume
that $P$ is a scheme by replacing the algebraic space $P$
by a presentation. As $X$ is locally of finite type,
we can assume that $P$ is in $\bf{lfschemes}/\C$, by
replacing $P$ by an open cover of $P$. We then 
take $Y^0\bu$ to be a resolution of 
singularities of $P$. We view $Y^0\bu$ as a 
$0$-truncated simplicial space.
 Note that a smooth morphism 
locally admits sections and a resolution
of singularities is proper and surjective so 
\[
 Y^0\bu{\rightarrow}\sX
\]
is a hypercover.

$k=1$. Let $Z_1 = (\cosk(Y^0\bu/\sX))_1$. We replace $Z_1$ 
by an open affine cover and then take a resolution of singularities
of this cover to obtain a smooth scheme $N_1$ in
$\bf{lfschemes}/\C$, and a map 
\[
 {\beta}:N_1{\rightarrow}Z_1.
\] 
Apply \ref{P:hyperconstruct} to the
triple $(Y^0\bu,N_1,{\beta})$ to obtain a 
smooth $1$-truncated split simplicial scheme
$Y^1\bu$.

Inductively one produces for each $k$ a split 
$k$-truncated simplicial scheme $Y^k\bu$ and
an augmentation $Y^k\bu{\rightarrow}\sX$ such that

\noindent
(1) The augmentation is a hypercover.

\noindent
(2) $Y_i^k$ is in $\bf{lfschemes}/\C$.

\noindent
(3) $Y_i^k$ is smooth over $\C$.

\noindent
(4) $\sk_{k-1}(Y^k\bu)=Y^{k-1}\bu$.

The last condition means that
\[
 Y_i^i = Y_i^{i+1}=\cdots
\]

We define $Y_i^\infty$ to be this stable value
of $Y_i^*$. The $Y_i^\infty$ fit together to
form a simplicial scheme that is in fact our required
hypercover
\[
 Y\bu=Y\bu^\infty{\rightarrow}\sX.
\]

A \emph{compactification} of a simplicial scheme $X\bu$
is a simplicial scheme $\bar{X}\bu$ and a morphism 
$j:X\bu{\hookrightarrow}\bar{X}\bu$ such that each of the
maps $j_n$ are compactifications.

A \emph{divisor} $D\bu$ on a smooth simplicial scheme $X\bu$
is a closed simplicial subscheme 
\[
 D\bu{\hookrightarrow}X\bu
\]
such that each of the morphisms $D_n{\hookrightarrow}X_n$
is a divisor. We say that $D\bu$ has \emph{simple
normal crossings} if each of the $D_n$ do.

\begin{theorem}
\label{T:hypercovers}
Let $\sX$ and $\sY$ be algebraic stacks locally of
finite type. 

\noindent
(i) We can construct a hypercover $X\bu\ri\sX$ 
with $X\bu$ smooth and a smooth compactification $\bar{X}\bu$ 
of $X\bu$ such that  $\bar{X}\bu\setminus X\bu$ is
a divisor with simple normal crossings and both
of these simplicial schemes are in $\bf{lfss}_\infty$.

\noindent
(ii) If we have two such hypercover-compactification pairs
$(X\bu,\bar{X}\bu)$ and $(X'\bu,\bar{X'}\bu)$ we can find a
third pair $(Z\bu,\bar{Z}\bu)$ that satisfies the conditions
of (i) and fits into a diagram
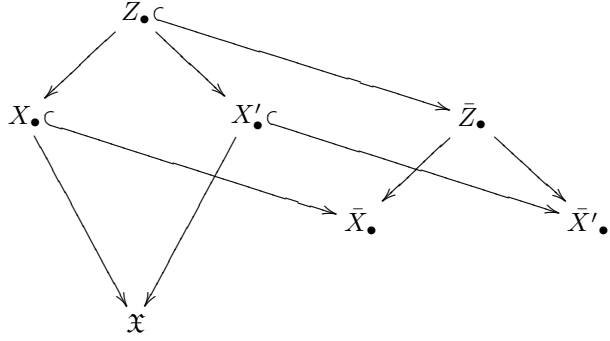
\[
\xymatrix{
     &  Z\bu  \ar@{^{(}->}[drrr] \ar[dr] \ar[dl] &        &      &     &   \\
X\bu \ar@{^{(}->}[drrr] \ar[ddr] &         &  X'\bu \ar[ddl]
     \ar@{^{(}->}[drrr] &            &   
\bar{Z}\bu \ar[dr]  \ar[dl]     &   \\
     &         &        & \bar{X}\bu &            &\bar{X'}\bu \\
     &   \sX   &        &            &            &           
}
\] 

\noindent
(iii) Let $F:\sX\ri\sY$ be a morphism. Then there exists 
hypercover-compactification pairs $(X\bu,\bar{X}\bu)$ 
and $(Y\bu,\bar{Y}\bu)$as in (i) for $\sX$ and $\sY$ 
respectively, along with morphisms
\[
X\bu\ri Y\bu\quad\bar{X}\bu\ri\bar{Y}\bu
\] 
and a commutative diagram
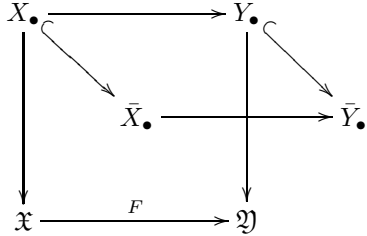
\[
\xymatrix{
X\bu \ar[dd] \ar[rr] \ar@{^{(}->}[dr]  &           &  Y\bu  \ar[dd]
\ar@{^{(}->}[dr] &   \\
       &\bar{X}\bu \ar[rr]&        & \bar{Y}\bu \\
\sX \ar[rr]^F   &           & \sY
}
\]

\end{theorem}

\begin{proof}
The proofs are analogous to those in
\cite{HodgeII}. For the convenience of the
reader we outline some of the proofs. 

(i) If $X$ is a scheme that is a disjoint union of smooth, separated,
finite type schemes
over $\C$, we may find a compactification of it
$\bar{X}$ by \cite{nagata:62}. We may assume that the compactification
is smooth and $\bar{X}\setminus X$ is a simple normal crossings 
divisor by \cite{hironaka:64}. The result will know follow from
the ideas in the discussion above.

(ii) The proof of this result is similar to that of (iii) so we only
give the proof of (iii).

(iii) Let $Y\ri\sY$ be a presentation of $\sY$. We may assume
that $Y$ is a disjoint union of separated schemes of finite type
over $\C$. The stack $\sX\times_{\sY}Y$ is algebraic and 
\[
\sX\times_{\sY}Y\ri\sX 
\]
is a representable surjective and smooth morphism. So a presentation
for this stack gives a presentation for $\sX$ by composition. So we
obtain a diagram
\[
\xymatrix{
X \ar[r]^f \ar[d]& Y \ar[d] \\
 \sX & \sY
}
\]
where the two vertical arrows are of universal cohomological descent and
$X$ and $Y$ are in $\bf{lfschemes}/\C$. We may further assume that
$X$ and $Y$ are smooth. To do this, first resolve $Y$ to 
$Y'$ and then resolve $X\times_{Y}Y'$ and note that the
projection $X\times_{Y}Y'\ri X$ is of universal cohomological descent.

We claim that there are smooth
compactifications of $X$ and $Y$ denoted $\bar{X}$ and $\bar{Y}$ 
respectively such that $f$ extends to a morphism
\[
\bar{f}:\bar{X}\ri\bar{Y}
\]
and 
\[
\bar{X}\setminus X \quad \bar{Y}\setminus Y
\]
are simple normal crossings divisors. To do this 
choose any compactifications $\bar{Y}$ of $Y$ and
$\bar{X'}$ of $X$. Let
\[
\bar{{\Gamma}_f}{\subseteq}\bar{X'}\times\bar{Y}
\]
be the closure of the graph of $f$. It is compact
and after applying \cite{hironaka:64} to it we may assume
that it is compact and the complement of the inclusion
$X{\subseteq}\bar{{\Gamma}_f}$  has simple normal crossings. 
We take $\bar{X}=\bar{{\Gamma}_f}$.
 This proves the claim.

We take $X_0=X$, $Y_0=Y$, $\bar{X}_0=\bar{X}$  and
$\bar{Y}_0=\bar{Y}$. To construct the next level of the required
simplicial
schemes form a diagram
\[
\xymatrix{
N' \ar[r]^{f_1} \ar[d]^{p'} & N \ar[d]^p \\
{\rm cosk}(X/\sX)_1 \ar[r] & {\rm cosk}(Y/\sY)_1,
}
\]
where $N$ and $N'$ are smooth schemes in $\bf{lfschemes}/\C$ and the vertical arrows are
of universal cohomological descent. We may compactify 
$N$ and $N'$ as above, so that $f_1$ extends to a morphism on the compactifications.
Now apply \ref{P:hyperconstruct} as in the discussion preceding
this theorem. One
continues
by induction and the required diagram is constructed.
\end{proof}

Consider the category whose objects are  pairs $(X\bu,\bar{X}\bu,)$ where
$X\bu$ and $\bar{X}\bu$ are smooth simplicial schemes
in ${\bf lfss}_\infty$ and ${\bar{X}}\bu$ is a compactification
of $X\bu$ with simple normal crossings on the boundary. 
We will now construct a functor from this category to
$\Q$-mixed Hodge structures. The underlying
vector space of this mixed Hodge structure will be
$\coh{*}{X\bu,\Q}$.

Once this functor is constructed, the above theorem will show
that a stack $\sX$ has a canonical functorial mixed Hodge 
structure. Note that a morphism of mixed Hodge structures
that is an isomorphism on underlying vector spaces is
in fact an isomorphism of mixed Hodge structures, so (ii)
above shows that the construction is independent of the
choice of hypercover-compactification. Functoriality follows
from (iii).

There is one \underline{very} minor complication here, as
$\coh{i}{X\bu,\Q}$ may not be of finite type we may not
directly directly apply \cite{HodgeII} and \cite{HodgeIII}.
However, we claim that once the definitions of these
papers are relaxed as outlined below the results of these
papers still hold.

An {\it infinite $\Q$-Hodge structure of weight $n$} is a $\Q$-vector
space $V$ and a finite decreasing filtration $F$ on
$V\otimes_\Q\C=V_\C$ such that the filtrations
$F$ and $\bar{F}$ are $n$-opposed, that is
\[ \gr{p}{F}\gr{q}{\bar{F}}(V_\C)= 0 \]
for $p+q\ne n$. We do not require that $V$ 
be finite dimensional.

\noindent
An {\it infinite $\Q$-mixed Hodge structure} consists of
the following data:
 
\noindent
(i) a $\Q$-module $V$.

\noindent
(ii) a finite increasing filtration $W$ on 
$V$, called the weight filtration.

\noindent
(iii) a finite decreasing filtration $F$ on
$V\otimes_\Q \C= V_\C$ called  the Hodge filtration.

This data is required to satisfy the following axiom:
$F$ induces a weight $n$ infinite Hodge structure 
on $\gr{W}{n}(V)$.

A morphism $f:V\rightarrow V'$ of infinite mixed Hodge
structures is a map of Abelian groups that induces
maps that are compatible with the filtrations.

\noindent
A weight $n$ {\it infinite Hodge complex} consists of 

\noindent
(${\alpha}$) A complex $K\ubu$ of $\Q$-modules.

\noindent
(${\beta}$) A filtered complex $(K_\C\ubu,F)$ in
$D^+F(\C)$ and an isomorphism
\[
K\ubu\otimes\C\isoarrow K_C\ubu\quad\text{in }
D^+(\C).
\] 

This data is required to satisfy the following axiom

\noindent
 For all $k$, the filtration on $\coh{k}{K\ubu_\C}$
induced by $F$,
defines a weight $n+k$ infinite Hodge structure.

In the above $D^+F(\C)$ is the filtered derived 
category as defined in \cite{HodgeIII}. In particular
the filtration $F$ is biregular, that is it a finite
filtration on each component of the complex
$K_\C\ubu$.

\noindent
An {\it infinite mixed Hodge complex} consists of

\noindent
($\alpha$) A filtered complex $(K,W)$ of $\Q$-vector
spaces in $D^+F(\Q)$.

\noindent
($\beta$) A bifiltered complex $(K_\C\ubu,W,F)$
a complex of $\C$ vector spaces, $W$ an increasing
biregular filtration, $F$ a decreasing biregular 
filtration and an isomorphism
\[
\C\otimes_\Q K\ubu\isoarrow K\ubu_\C\quad
\text{ in } D^+F(\C).
\] 
This data is required to satisfy the following axiom:

\noindent
The data consisting of the complex $\gr{W}{n}{K\ubu_\Q}$ and
the quasi isomorphism
\[
\gr{W}{n}{K\ubu}\otimes\C\isoarrow \gr{W}{n}{K_\C\ubu}
\]
is a weight $n$ infinite Hodge complex.

We will now proceed to show that the cohomology of an infinite
mixed Hodge complex inherits a canonical infinite mixed Hodge
structure. We first need to recall some facts from
\cite{HodgeII}.

Let $(K\ubu,W,F)$ be a bifiltered complex. On the 
terms $E_r^{pq}(K\ubu,W)$ of the spectral sequence
associated to the filtered complex $(K\ubu,W)$, we have
three filtrations induced by $F$:

(i) {\it The first direct filtration} $F_d$, is formed
by viewing $E^{pq}_r$ as a quotient of a subobject of
$K^{p+q}$. 

(ii) {\it The second direct filtration} $F_{d^*}$, is formed
by viewing $E^{pq}_r$ as a subobject of a quotient object
of $K^{p+q}$. 

(iii) {\it The recursive filtration} $F_r$, is formed by
defining 
\begin{eqnarray*}
E_0^{pq} & : & F_r=F_d=F_{d*} \mbox{ (See below)}\\
E_r^{pq} & : & \mbox{ The filtration induced by the direct filtration
 on }E_{r-1}^{pq} 
\end{eqnarray*}

\begin{prop}
\label{P:3filt}
We have:

(i) On $E_0$ and $E_1$ the three filtrations coincide.

(ii) The differentials $d_r$ are compatible with $F_d$ and 
$F_{d^*}$. 

(iii) We have $F_d{\subseteq}F_r{\subseteq}F_{d^*}$.
\end{prop}

\begin{proof}
See \cite[pg. 17]{HodgeII}.
\end{proof}

\begin{theorem}
\label{T:3filt}
Let $(K\ubu,W,F)$ be a bifiltered complex. We let
$E_r^{pq}=E_r^{pq}(K\ubu,W)$ be the terms of the spectral
sequence. Suppose that $F$ is biregular and 
for $0\le r \le r_0$ the differentials $d_r$ are
strictly compatible with $F_r$
Then on $E_{r_0+1}$ we have $F_d=F_r=F_{d*}$.
\end{theorem}

\begin{proof}
See \cite[pg. 18]{HodgeII}.
\end{proof}

Given a complex $K\ubu$ with an increasing filtration 
$W$. We define a new shifted filtration ${\rm Dec}W$
on $K\ubu$ by
\[
{\rm Dec}W_n K^i=W_{n-i} K^i.
\]

\begin{theorem}
Let $(K\ubu,W,{\alpha},K\ubu_\C,F)$ be an infinite
mixed Hodge complex. Then
${\rm Dec}(W)$ and $F$ induce a mixed Hodge structure
on $\coh{i}{K\ubu}$
\end{theorem}

\begin{proof}
Consider the decreasing filtration
$\tilde{W}$ on $K$ defined by 
\[
\tilde{W}^p = W_{-p}.
\]
This filtration gives a spectral sequence 
with
\[
E_1^{pq} = \coh{p+q}{\gr{W}{-p}(K)}
\]
abutting to $\coh{p+q}{K}$. By 
Proposition \ref{P:3filt} the three filtrations
on $E_1^{pq}$ coincide and the differential is
compatible with this filtration. As $d_1$ is defined
over $\Q$, this differential is compatible with
the conjugate filtration and therefore is strictly compatible
with the filtration. So 
\[
d_1:E_1^{pq}\ri E_1^{p+1,q}
\]
is a morphism of Hodge structures of weight $q$. 

Hence $E_2^{pq}$ has a weight $q$ Hodge structure. By 
Theorem \ref{T:3filt} the three filtrations coincide
on $E_2$ and $d_2$ is compatible with it. As before we
conclude that $d_2$ is strictly compatible with
this filtration. However
\[
d_2:E_2^{pq}\ri E_2^{p+2,q-1}
\]
is a morphism of Hodge structure of different weights
so it vanishes. Hence $E_2^{pq} = E_\infty^{pq}$
and so
\[
\gr{W}{-p}\coh{p+q}{K}
\]
has a weight $q$ Hodge structure. One checks that
\[
\gr{{\rm Dec}}{q}\coh{p+q}{K}=\gr{W}{-p}\coh{p+q}{K}
\]
and we are done.
\end{proof}
 
One can now proceed to define infinite complexes of sheaves
as in \cite[pg. 28-38]{HodgeIII}. The results will carry over
verbatim to this setting. In particular, the analogue of 
Proposition (8.1.20)
constructs a functorial mixed Hodge structure on the
cohomology of a hypercover-compactification pair.

\section{The Cohomology of the Ind Scheme of Matrix Divisors}

For the remainder of this paper $X$ is a smooth geometrically
connected projective curve defined over our ground field $k$.

The primary purpose of this section  is to
recall the results in \cite{bifet:94} regarding the cohomology
of $\bdiv$ and provide a bound on the codimension of the
complement $\bdiv^{ss}\setminus\bdiv^s$. 

Let ${\Lambda}$ be the partially ordered set of effective divisors on 
$X$. Fix $D\in{\Lambda}$ and consider the functor
\begin{eqnarray*}
\div^{r,d}(D)^\flat : \schemes/k{\rightarrow}{\bf sets}
\end{eqnarray*}
whose $S$-points are equivalence classes of 
inclusions
\[
 \F{\hookrightarrow}{\mathcal O}_{X\times S}(D)^r,
\]
where $\F$ is a family of rank $r$ degree $d$ bundles on
$X\times S$. This functor is representable by a Quot
scheme that we denote by
\[
 \div^{r,d}(D)=\div(D).
\]
These Quot schemes fit together to form an ind scheme 
denoted
\[
 \bdiv^{r,d}=\bdiv.
\]

Let ${\bf m}=(m_1,m_2,\ldots,m_r)$ be a partition of the 
integer $r.\deg D - n$, by non negative integers. 
Then the product of Hilbert schemes
of points 
\[
 H^{\bf m}=\hilb^{m_1}(C)\times \hilb^{m_2}(C)\times \ldots 
 \hilb^{m_r}(C)
\]
sits canonically inside of $\div(D)$. Recall that over 
an algebraically closed field, the Hilbert scheme of
points of a smooth curve is just a symmetric power
of the curve.

The torus ${\mathbb G}_m^r$ acts on $\div(D)$ and the
above products of Hilbert schemes are clearly fixed 
by this action. The converse is also true. 

\begin{theorem}
(i) The fixed points of this action are precisely the
schemes $H^{\bf m}$ as ${\bf m}$ varies over all
partitions of $r.\deg D -n$. 

\noindent
(ii) The cohomology of $\bdiv$ stabilizes and its Poincare
polynomial is given by
\[
 P(\bdiv;t) = \frac{\prod_{i=1}^r (1+t^{2j-1})^{2g}}
{(1-t^{2r})\prod_{i=1}^{r-1}(1-t^{2j})^2}.
\]
The fact that the cohomology stabilizes means that 
the inverse limit
\[
 \varprojlim_{\Lambda}\coh{i}{\div(D),\Q}
\]
is in fact finite. 

\noindent
(iii) When $k=\C$ the Hodge-Poincare polynomial of
$\bdiv$ is 
\[
 P_H(\bdiv;x,y) = \frac{(1+x)^g(1+y)^g}{(1-x^ry^r)}
\prod_{i=1}^{r-1}\frac{(1+x^{i+1}y)^g(1+xy^{i+1})^g}
{(1-x^iy^i)^2}.
\]
\end{theorem}

\begin{proof}
The first part is proved in 
\cite{bifet:89}. The second part follows from the
first by some  theorems of  A. Bialynicki-Birula and
some deformation theory.
For details see \cite{bb:73},\cite{bb:74} and \cite[Proposition 4.2]{bifet:94}.
The last part follows by noting that the
Bialynicki-Birula decomposition is compatible with, 
amongst other things, Hodge theory. A nice exposition
of these ideas can be found in \cite{delbano:01}. The formula
we have written down follows directly from Proposition 4.4
of this paper.
\end{proof}

For a vector bundle $\E$ on $X$ with rank $r$ and
degree $d$. Its Harder-Narasimhan
\[
 \E_1{\subseteq}\E_2{\subseteq}\ldots{\subseteq}\E_l=E.
\]
filtration is unique. So the sequence of pairs of numbers 
$(r_1,d_1),(r_2,d_2)\ldots,(r_l,d_l)$, where $r_i$ 
is rank of $E_i$ and $d_i$ its degree, is unique. If these
points are plotted in ${\mathbb R}^2$ and the line 
segments $(r_i,d_i)$ to $(r_{i+1},d_{i+1})$ joined one obtains
a polygonal curve from the origin to 
$(r,d)$ such that the slope of each successive line segment 
decreases. Such a curve will be called a \emph{Shatz polygon}
for $(r,d)$. We denote the set of Shatz polygons for
$(r,d)$ by ${\mathcal P}^{r,d} = {\mathcal P}$. 
If one thinks of these polygons as graphs of
functions $[0,r]{\rightarrow}{\mathbb R}$ then this 
collection has a natural partial order determined by the
partial order on the set of functions with domain $[0,r]$
and codomain ${\mathbb R}$. For a vector bundle
$\E$, we let $s(\E)$ denote  its Shatz polygon.

Now consider a family of vector bundles $\E$ on $X\times T$ 
of rank $r$ and degree $d$, with $T$ in 
${\bf lfschemes}/k$. 
Fix a Shatz polygon $P$ for $(r,d)$ and
recall the following results:
\noindent
(i) The locus 
\[
 T^P=\{t\in T| s(\E_t)>P \}
\]
is closed.

\noindent
(ii) The locus
\[
 \{t\in T| s(\E_t)=P \}
\]
is closed in the open set
$T\setminus T^P$. 
To prove these statements one considers the relative flag scheme
${\rm FLag}^P(\E/T)$ over $T$, whose fiber over $t\in T$ 
is a parameter space for flags of $\E_t$ with rank and degree
data specified by $P$. It is proper over $T$ so it
has closed image in $T$. The above results follow by use of this
fact. Complete details can be found in \cite{bruguieres:83}.

Denote by $\div^P(D)$ the open locus inside 
$\div(D)$ parameterizing subbundles of 
${\mathcal O}_X(D)^r$ whose Shatz polygon is not
bigger than $P$, i.e the complement of the
closed set in (i) defined by taking 
$T=\div(D)$. We can consider the corresponding 
ind schemes $\bdiv^P$. We denote by 
$\bdiv^{ss}$ the semistable locus, corresponding
to taking $P$ equal to the straight line from
$(0,0)$ to $(r,d)$.

Denote by $S^P(D)$ the locally closed locus inside
$\div(D)$ parameterizing bundles with Shatz polygon exactly $P$. 
These fit together to form an ind scheme
${\bf S}^P$. For $\deg D$ large enough $S_P(D)$ is
smooth. If $P$ has vertices
\[
 (r_0=0,d_0=0),(r_1,d_1),\ldots,(r_l=r,d_l=d)
\]
and $\deg D$ large
then the codimension of this stratum is given by
\[
 d_P = \sum_{i<j} r_ir_j({\mu}_i-{\mu}_j + g-1)
\]
where ${\mu}_i = d_i/r_i$.

\begin{theorem}
Let $P$ be a Shatz polygon with vertices
\[
 (r_0=0,d_0=0),(r_1,d_1),\ldots,(r_l=r,d_l=d).
\]
Set $r_i'=r_i-r_{i-1}$ and $d_i'=d_i-d_{i-1}$. There
is a closed immersion
\begin{eqnarray*}
 {\delta}& :& \bdiv^{r_1',d_1',ss}\times
            \bdiv^{r_2',d_2',ss}\times\ldots\times
            \bdiv^{r_l',d_l',ss}\rightarrow {\bf S}^P \\
 & & (\E_1,\E_2,\ldots,\E_l)\mapsto \E_1\oplus \E_2 \oplus \E_l
\end{eqnarray*}
that induces an isomorphism in cohomology.
\end{theorem}

\begin{proof}
This is \cite[Proposition 7.1]{bifet:94}.
\end{proof}

Let $I$ be a subset of the collection of all matrix
divisors. We say that $I$ is \emph{open} if
$P\in I$ and $P'\le P$ implies $P'\in I$. If $P$ is a 
minimal element of the complement of $I$ then
$J=I\cup\{ P \}$ is also open.
If $I$ is open then the locus
\[
S^I = \cup_{P\in I} S^P
\]
is an open subset of $\div(D)$. 

\begin{theorem}
Suppose $P$ is a minimal element of the
complement of $J$ with $J$ open. Set
$I = J\cup \{P\}$.  
The Gysin sequences
\[
\ri\coh{i-2d_P}{{\bf S}^P,\Q}\ri \coh{i}{{\bf S}^I,\Q}\ri
\coh{i}{\bf{S}^J,\Q}\ri,
\]
split into short exact sequences. Hence
the following relation among Poincare
polynomials holds:
\[
 P(\bdiv;t)=\sum_{P\in{\mathcal P}} P({\bf S}^P;t)t^{2d_P}.
\]
\end{theorem}

\begin{proof}
See \cite[Proposition 10.1]{bifet:94}.
\end{proof}

The above three theorems yield recursive 
formulas for the Hodge and Betti numbers of
the ind varieties of matrix divisors associated
to Shatz polygons.

In the remainder of this section we provide a dimension bound for
the complement $\bdiv^{ss}\setminus\bdiv^s$.

We consider pairs of sequences of integers
\[
(\ur,\ud)=((r_1,r_2,\ldots,r_l),(d_1,d_2,\ldots,d_l))
\]
satisfying the following conditions

\begin{equation}
\label{e:rel}
0<r_1<r_2<\cdots< r_l=r\qquad d_i=\frac{dr_i}{r}.\tag{$\sharp$}
\end{equation}
We denote by $\flag^\prd(D)$ the scheme representing the
functor
\[
T\mapsto\{\E_1{\subseteq}\E_2{\subseteq}\ldots{\subseteq}\E_l 
{\subseteq}\structd |
rk\E_i=r_i\ deg\E_i=d_i\}
\]

See \cite{bifet:94} for the existence of such a scheme.
There is a proper morphism 
\[
\pi^\prd:\flag^\prd(D)\ri \divf(D).
\]
There is an open subset 
\[
JH^\prd(D){\subseteq}\flag^\prd(D)
\]
parameterizing semistable flags with $\E_i/\E_{i-1}$ 
a stable bundle for all $i$. By the 
existence of Jordan-Holder filtrations we have
\[
\divf^s(D)=\div(D)\setminus \cup_{\prd}
\pi^\prd(JH^\prd).
\]
To find a dimension bound on the complement
\[
\div(D)\setminus \divf^s(D)
\]
we need only bound the dimensions of each
of the open sets $JH^\prd(D)$. We will
show
\begin{theorem}
\[
{\rm dim}JH^\prd(D)\le r^2{\rm deg}D-rd-(g-1)(r-1),
\]
where $g$ is the genus of the curve.
\end{theorem}

\begin{proof}
Consider a point
\[
E_1{\subseteq}E_2{\subseteq}\cdots{\subseteq}E_l{\subseteq}\structd
\]
of $JH^\prd(D)$. Following \cite{bifet:94} we denote by 
$\tE_i$ the sheaf $\structd/\E_i$.
From \cite{bifet:94}, the  tangent space to
$JH^\prd(D)$ at 
the above point is identified with the vector subspace of
\[
\hom{}{\E_1,\tE_1}\oplus\hom{}{\E_2,\tE_2}\oplus\cdots
\oplus\hom{}{\E_l,\tE_l}
\]
consisting of $l$-tuples $(x_1,x_2,\ldots,x_l)$ satisfying the
following condition:

\noindent
The image of $x_i$ and $x_{i+1}$ agree in 
$\hom{}{\E_i,\tE_{i+1}}$. (See\cite{bifet:94}.)

We have exact sequences 
\begin{eqnarray*}
0\ri \E_1\ri R_{i+1}\ri L_i\ri 0 \\
0\ri L_i\ri \tE_i\ri \tE_{i+1}\ri 0
\end{eqnarray*}
where $L_i$ is a stable bundle of rank
$r_{i+1}-r_i$. These sequences give rise to
long exact sequences
\[
0\ri\hom{}{\E_i,L_i}\ri \hom{}{\E_i,\tE_i}\ri
\hom{}{E_i,\tE_{i+1}}\ri\cdots
\]
and
\[
0\ri\hom{}{L_i,\tE_{i+1}}\ri\hom{}{\E_{i+1},\tE_{i+1}}
\ri\ext{1}{}{L_i,\tE_{i+1}}\ri\cdots
\]
As $L_i$ is stable and $E_i$ is semistable
$\hom{}{E_i,L_i}=0$ and for $\deg D$
large enough $\ext{1}{}{L_i,\tE_{i+1}}$ vanishes.
It follows that 
\begin{eqnarray*}
{\rm dim}JH^\prd(D)&\le&
{\rm dim}\hom{}{\E_1,\tE_1} + {\rm dim}\hom{}{L_1,\tE_2} +
{\rm dim}\hom{}{L_2,\tE_3} \\ & &  + \ldots + 
 {\rm dim}\hom{}{L_{l-1}\tE_l}.
\end{eqnarray*}
In bounding the right hand side above, we will freely make use
of the formulas, theorems and notations of \S5 and \S15 of
\cite{fulton:98}. We have
\begin{eqnarray*}
Td(C)& =& 1 + \frac{1}{2}c_1(-K) \\ 
ch(\E_1^\dual)& =& r_1 - c_1(E_1) \\
ch(\tE_1)& =& r-r_1+rc_1(D)-c_1(E_1) \\
ch(\tE_1\otimes \E_1^\dual)& =& r_1(r-r_1) + r_1rc_1(D) - rc_1(D).
\end{eqnarray*}
Hence,
\[
{\chi}(\tE_1\otimes E_1^\dual) = 
r_1(r-r_1)(1-g) + r_1r\deg D -r_1d.
\]
Similarly,
\begin{eqnarray*}
ch(L_i^\dual\otimes\tE_{i+1}) &=&
(r-r_{i+1})(r_{i+1}-r_i) + (r-r_{i+1})c_i(\E_i)
(r_i-r)c_1(\E_{i+1})\\ && +(r_{i+1}-r_i)rc_1(D).
\end{eqnarray*}
Hence, using \ref{e:rel},
\begin{eqnarray*}
\chi(L_i^\dual\otimes\tE_{i+1}) &=& 
(r-r_{i+1})(r_{i+1}-r_i)(1-g) + \\
& & (r-r_{i+1})d_i + 
(r_i-r)d_{i+1} + (r_{i+1}-r_i)r\deg D \\
& = & (r_{i+1}-r_i)r\deg D + (r-r_{i+1})(r_{i+1}-r_i)(1-g) \\
& & + (r_i - r_{i+1})d.
\end{eqnarray*}
So 
\begin{eqnarray*}
{\rm dim} JH^\prd(D)& \le & r_1(r-r_1)(1-g) + r_1r\deg D -r_1d + \\
   & & (r_2-r_1)r\deg D + ((r-r_2)(r_2-r_1)(1-g) + (r_1-r_2)d + \\
   & & (r_3-r_2)r\deg D + (r-r_3)(r_3-r_2)(1-g) + (r_2 -r_3)d + \\
   & & \cdots + \\
   & & (r_l - r_{l-1})r\deg D + (r-r_l)(r_{l}-r_{l-1})(1-g) + \\
   & & (r_{l-1}-r_l)d \\
   & = & r^2\deg D -rd + \\
   & & (1-g)(r_1(r-r_1)+(r-r_2)(r_2-r_1) + (r-r_l)(r_l-r_{l-1}) \\
   &\le & r^2\deg D -rd + (1-g)(r-1).
\end{eqnarray*}
To see why the last inequality holds first observe that for 
integers $s$ and $t$ with $1\le s<t$ we have
\[
\frac{s(t-s)}{t-1}\ge 1.
\]
Since $1-g<0$, the inequality in question is equivalent to
 showing that
\[
r_1(r_2-r_1)+ r_2(r_3-r_2) + \ldots + r_{l-1}(r-r_{l-1})\ge r-1
\]
which follows from the above inequality.
\end{proof}

\begin{coro}
\label{C:stable}
The inclusion 
\[
\divf^s(D)\hookrightarrow \div(D)
\]
induces an isomorphism in cohomology
\[
\coh{i}{\div(D),\Q}\isoarrow 
\coh{i}{\divf^s(D),\Q}
\]
for $i < 2(g-1)(r-1)$.
\end{coro}

\begin{proof}
We calculate the dimension of 
$\divf(D)$ to be $r^2\deg D -rd$. 
The result follows from the Gysin sequence
and the dimension bound above, see for
example \cite[pg. 268]{milne:80}.
\end{proof} 

\section{The Cohomology of the Stack}

Let $\E$ be a family of vector bundles on $X\times S$
with $S/k$ a smooth scheme. We say that the family
$\E$ is \emph{complete} if the Kodaira-Spencer map
\[
 T_sS{\rightarrow}\ext{1}{}{\E_s,\E_s}
\]
is surjective.

\begin{lemma}
Let $S$ and $T$ be schemes smooth over $k$ and let
$\E_S$ (resp. $\E_T$) be a complete family of bundles on
$X\times S$ (resp. $X\times T$). Assume also that the
induced maps $S{\rightarrow}\bun$ and $T{\rightarrow}\bun$
are smooth. Then the induced family on
$S\times_{\bun}T$ is complete.
\end{lemma}

\begin{proof}
This is mostly a matter of unwinding definitions. Recall
that if $A$ is a $k$-algebra then an $A$-point on
$S\times_{\bun}T$ consists of a triple $(s,t,{\alpha})$
where $s$ (resp. $t$) is an $A$-point of $S$ (resp. $T$)
and ${\alpha}$ is an isomorphism
\[
 {\alpha}:(s\times 1)^*\E_S\stackrel{\sim}{\rightarrow}
          (t\times 1)^*E_T.
\] 

So consider a closed point $(s_0,t_0,{\alpha}_0)$ of
the fibered product. We have a diagram of
Kodaira-Spencer maps
\[
 \xymatrix{
T_{s_0}S \ar[r] & \ext{1}{}{\E_{s_0},\E_{s_0}} \ar[d]^{\sim} \\
 T_{t_0}T \ar[r] & \ext{1}{}{\E_{t_0},\E_{t_0}},
}
\]
where the vertical arrow is an isomorphism and the horizontal
maps are surjective. Fixing an extension class  and
chooses $k[{\epsilon}]$-points of $S$ and $T$ lying above it.
Call these point $s$ and $t$ respectively. The bundles
\[
 (s\times 1)^*\E_S\quad\text{and}\quad(t\times 1)^*E_T
\]
are isomorphic as they correspond to the same extension 
class. It is possible to choose an isomorphism between 
these bundles that
restricts to ${\alpha}$ upon specialization to the
closed point of $k[{\epsilon}]$. Such an isomorphism
gives a $k[{\epsilon}]$-point of the fibered product
that maps onto the extension class we chose earlier.
\end{proof}

We recall how a presentation of $\bun$ was constructed
in \cite{laumon:00}. 
Let $p(x)=rx + d +r(1-g)$. For every
integer $m$ we define an open subscheme
\[
 Q^m{\hookrightarrow}{\rm Quot}
({\mathcal O}_X^{p(m)},p(x+m))
\]
by requiring that

\noindent
(i) The quotients parameterized by $Q^m$ be vector bundles.

\noindent
(ii) For every $T$-point of $Q^m$ defined by the
family
\[
 {\tau}:{\mathcal O}_{X\times T}^{p(m)}{\rightarrow}\F
{\rightarrow}0,
\]
we have $R^1{\pi}_{T,*}\F=0$ and
\[
 \pi_{T,*}:{\mathcal O}_{X\times T}^{p(m)}
\stackrel{\sim}{\rightarrow}\pi_{T,*}\F
\]
is an isomorphism.

It follows from (ii) that if the
quotient
\[
 0{\rightarrow}\G{\rightarrow}
{\mathcal O}_X^{p(m)}{\rightarrow}\F
{\rightarrow}0
\]
represents a point of $Q^m$ then we have
$\coh{1}{\F\otimes \G\dual}=0$, i.e $Q^m$ is
smooth. 

We have maps
\begin{eqnarray*}
 Q^n{\rightarrow}\bun \\
 \F \mapsto \F(-n).
\end{eqnarray*}
Then 
\[
 Q = \coprod_m Q^m{\rightarrow}\bun
\]
is a smooth presentation.

\begin{proposition}
\label{P:complete}
The family on $Q$ is complete and hence,
by the lemma,
\[
 \cosk(Q/\bun)
\]
is a simplicial algebraic space each of whose
components defines a smooth family. 
\end{proposition}

\begin{proof}
Let 
\[
0{\rightarrow}\G_n{\rightarrow}
{\mathcal O}_{X\times Q^n}^{p(n)}{\rightarrow}
\F_n{\rightarrow}0,
\]
be the universal family on $Q_n\times X$. 
The Kodaira-Spencer map is identified with
the connecting homomorphism
\[
 \hom{}{\G_n,\F_n}{\rightarrow}\ext{1}{}{\F_n,\F_n}=
\ext{1}{}{\F_n(-n),\F_n(-n)}.
\]
The next term in the sequence vanishes and the result
follows.
\end{proof}

We note the following theorem from \cite[pg. 206]{lepotier:97}.

\begin{theorem}
Let $E$ be a complete family of vector bundles $X$ parameterized
by $S$. Assume $S$ is smooth. Let
\[
P=((d_1,r_1),\ldots,(d_l,r_l))
\]
be a Shatz polygon for $(r,d)$. 
Then the subvariety $S_P$ of $S$ parameterizing 
bundles with polygon $P$ is locally closed and has codimension
\[
cod(P)\stackrel{\text{defn}}{=}
\sum_{i<j} r_ir_j({\mu}_i-{\mu}_j + g - 1),
\]
where $\mu_i = d_i/r_i$.
\end{theorem}

Let $\bun^P$ be the open substack of $\bun$ 
parameterizing bundles whose Shatz polygon is not
bigger than $P$.

\begin{theorem}
\label{T:stackshatz}
We have
\[
 \varprojlim_P\coh{i}{\bun^P,\Q} = \coh{i}{\bun,\Q},
\]
and in fact the limit on the left stabilizes.
\end{theorem}

\begin{proof}
First some notation, if $\E$ is a family
of rank $r$ degree $d$ bundles on $T\times X$,
denote by $T^P$ the open locus consisting
of points $t\in T$ such that $s(\E_t)$ is not
bigger than $P$. 
From definitions we have 
\[
 (S\times_{\bun} T)^P = (S^P \times_{\bun^P} T^P).
\]

For each fixed integer $i$ there are only finitely 
many Shatz polygons $Q$ with $cod(Q)<i$. 
Let $P_0$ be a Shatz polygon greater than all of
the Shatz polygons in this finite set. Let
$P\ge P_0$. It suffices to show that the natural map
\[
 \bun^P{\rightarrow}\bun
\]
induces an isomorphism on degree $i$ cohomology.
By \ref{P:ss} and \ref{P:complete}, it suffices to show
that if $\E$ is a complete family of vector bundles
on $X\times T$ then the natural inclusion
\[
 T^P{\hookrightarrow} T
\]
induces an isomorphism in cohomology of 
degree $j$ for all $j\le i$. But this follows by
the Gysin sequence and choice of $P$.
\end{proof}

The virtue of the above theorem is that the
family of bundles parameterized by $\bun^P$ 
is bounded. To see this last statement, note that
only finitely many Shatz polygons appear in $\bun^P$ 
and that the collection of bundles with a particular
Shatz polygon is bounded. We will now proceed to exploit this.

\begin{theorem}
\label{T:quasi}
The natural map
\[
 \bdiv{\rightarrow}\bun
\]
is a quasi-isomorphism, i.e it induces an 
isomorphism on cohomology groups.
\end{theorem}

The proof will rely on the following:
(See \cite[Lemma 8.2]{bifet:94} for a proof)

\begin{proposition}
Let $\E$ and $\F$ be rank $r$  bundles on
$X$ such that $\ext{1}{}{\E,\F}=0$. Then for 
any effective divisor $D$, the codimension $c_D$ of
the closed locus in $\hom{}{\E,\F(D)}$ consisting
of non-injective homomorphisms satisfies
\[
 c_D\ge\deg D.
\]
\end{proposition}

\begin{proof}(of the theorem) By \ref{T:stackshatz},
it suffices to show that the natural map
\[
 \bdiv^P{\rightarrow}\bun^P, 
\]
is a quasi-isomorphism. As this last stack is
of finite type, it suffices to show, by \ref{P:ss},
that for all schemes $T$ of finite type and all maps
$T{\rightarrow}\bun^P$ the map 
$p_T$ below is a quasi-isomorphism:
\[
 \xymatrix{
T\times_{\bun^P}\bdiv^P \ar[d] \ar[r]^(.7){p_T} & T \ar[d] \\
\bdiv^P \ar[r] & \bun^P. 
}
\]
Let $\F$ be the family of bundles on 
$X\times T$ defining the map 
$T{\rightarrow}\bun^P$. For $D$ large enough
we have $\coh{1}{\F_t\dual(D)}=0$ for each
$t\in T$. So, by the standard results on base change,
for $D$ large enough, we have that an  $S$-point of 
$T\times_{\bun^P}\div^P(D)$ consists of
a map ${\phi}:S{\rightarrow}T$ and an injection
\[
 ({\phi}\times 1)^*\F{\hookrightarrow}{\mathcal O}_{X\times S}(D)^r.
\]
Hence $T\times_{\bun^P}\div^P(D)$ is an open
subset of the vector bundle 
\[
 \pi_{T,*}({\mathcal Hom}(\F,{\mathcal O}_{X\times T}(D)^r.
\]
The result follows by the above lemma and
a Gysin sequence.
\end{proof}

\begin{corollary}
When $k=\C$ or a finite field, the cohomology 
of $\bun$ is pure and of the correct weight.
\end{corollary}

\section{The Cohomology of the Stack Versus that of the
Moduli Space}

\begin{proposition}
Let $G$ be a geometrically connected group scheme over $k$ and let
$f:P{\rightarrow}Y$ be a $G$-torsor. Then the local systems
\[
 Rf_{*} \Q\quad\text{ or for the etale site}\quad
Rf_*\Q_l,
\]
are constant.
\end{proposition}

\begin{proof}
This result is from
\cite[\S 1.4]{behrend:93}. For the convenience of
the reader give an outline of the ideas.

Firstly the action of $G$ on itself by left 
multiplication induces an action of $G$ on
$\coh{i}{G}$. This action is trivial as it comes
from an action of $G$ on the discrete spaces
\begin{eqnarray*}
 \coh{i}{G,\Z} \quad \text{if } k=\C \\
 \text{or for general }k \qquad
 {\rm H}^i_{et}(G\otimes_k \bar{k},\Z/l^n\Z).
\end{eqnarray*}
The fibration 
\[
g: P\times_G\coh{i}{G}{\rightarrow}Y
\] is hence
trivial over $Y$ and one shows that
$Rf_*\Q = Rg_*\Q$.
\end{proof}

\begin{proposition}
\label{P:glpgl}
Consider natural map
\[
 {\Phi}:\gl_n{\rightarrow}\pgl_n\times{\mathbb G}_m
\]
that is the projection on the first factor and the determinant
on the second factor. Then $\Phi$ is a quasi-isomorphism.
\end{proposition}

\begin{proof}
Consider the commutative diagram
\[
 \xymatrix{
 \gl_n \ar[r]^(.35)\Phi \ar[dr]^f & \pgl_n\times{\mathbb G}_m \ar[d]^g\\
  & \pgl_n,
}
\]
where $f$ and $g$ are the projections. It suffices to show
that map induced by $\Phi$ on the Leray spectral sequences
for $f$ and $g$ is an isomorphism at the $E_2$ level. We have
\[
 {\Phi}^*:E_2^{pq}(g)=\coh{p}{\pgl_n,\coh{q}{{\mathbb G}_m,\Q}} 
{\rightarrow}E_2^{pq}(f)=\coh{p}{\pgl_n,Rf_*^q\Q}
\]
By the above Proposition the local system on the right is
constant, and it suffices to observe that the $n$th power
map
\[
 {\mathbb G}_m{\rightarrow}{\mathbb G}_m
\]
induces an isomorphism in cohomology.
\end{proof}

Let $G$ be an algebraic group over $k$ acting on a
scheme $X$ over $\C$. Let the action be given by
\[
{\sigma}:X\times G\ri G.
\] 
The map $X{\rightarrow}[X/G]$ is a presentation
and we wish to describe the simplicial space
$\cosk(X/[X/G])$.
The $n$th  term of the simplicial scheme 
\[{\rm cosk}(X/[X/G])\] is of the form
\[
X\times\underbrace{G\times G\times\cdots\times G}_{n\text{ times}}
\]
The  $i$th face map is given by
\[
{\delta}_i(x,g_1,g_2,\ldots,g_n)=(x,g_1,\ldots,\hat{g_i},\ldots,g_n)\quad
\text{for }i>0,
\]
and
\[
{\delta}_i(x,g_1,g_2,\ldots,g_n)=(xg_1,g_1^{-1}g_2,\ldots,g_1^{-1}g_n)\quad
\text{for }i=0.
\]

Let $\bun^s$ be the open substack of $\bun$ parameterizing 
stable bundles. The following is well known:

\begin{prop}
There is a scheme $Q$ and
a commutative diagram of stacks
\[
\xymatrix{
[Q/\gl] \ar[r] \ar[d] & [Q/\pgl] \ar[d]\\
 \bun^s \ar[r] & \mo
}
\]
in which the two vertical arrows are isomorphisms.
\end{prop}

\begin{proof}
The $Q$ in the above theorem is the open
locus inside the quot scheme that is both
a presentation for $\bun^s$ and the GIT quotient
of it by $\pgl$ is the moduli space. For details,
see \cite{gomez:99}, Proposition 3.3.
\end{proof}

\begin{theorem}
There is an isomorphism
\[
H^*(\bun^s,\Q) \isoarrow H^*(\mo,\Q)\otimes H^*(B\gm,\Q).
\]
\end{theorem}

\begin{proof}
 We have
a map
\[
{\rm cosk}(Q/[Q/\gl]) \rightarrow 
{\rm cosk}(Q/[Q/\pgl])
\]
We define a map
\[
{\rm cosk}(Q/[Q/\gl]) \rightarrow {\rm cosk}(\text{point}/B\gm),
\]
 by  projecting onto
${\rm cosk}(\text{point}/B{\gl})$ and then taking the
determinant. We hence have an a map
\[
{\rm cosk}(Q/[Q/\gl]) \rightarrow 
{\rm cosk}(Q/[Q/\pgl]) \times 
{\rm cosk}(\text{point}/B\gm).
\]
We see that this map induces an isomorphism in cohomology 
by using the
standard spectral sequence \ref{P:ss}, and 
\ref{P:glpgl}.
\end{proof}

\begin{corollary}
The natural map 
\[
 \div^s{\rightarrow}\mo
\]
is a quasi-isomorphism. For $i<2(r-1)(g-1)$
the Betti and Hodge numbers of $\coh{i}{\mo}$ can be 
calculated. If $k$ is a finite
field or $\C$ these cohomology groups are pure and of
the correct weight.
\end{corollary}

\begin{proof}
Firstly, the natural map 
\[
 \bun^s{\rightarrow}\mo
\]
is a quasi-isomorphism, as the proof
of \ref{T:quasi} carries over to this 
case verbatim. The result now follows 
from the above theorem and \ref{C:stable}.
\end{proof}

\bibliographystyle{alpha}
\bibliography{./../latestbib/mybib,./../latestbib/sga}

\end{document}